\documentclass[12pt]{article}
\usepackage{latexsym, amssymb, amsmath, amscd, amsfonts, epsfig, graphicx, colordvi,verbatim,ifpdf,extarrows}
\usepackage{amsfonts, amsmath, amssymb}
\usepackage{amssymb,amsfonts,amsmath,latexsym,epsfig,cite, psfrag,eepic,color}
\usepackage{amscd,graphics}
\usepackage{latexsym, amssymb,  amsmath,amscd, amsfonts, epsfig, graphicx, colordvi,amsthm}

\usepackage{graphicx}
\usepackage{epstopdf}
\usepackage{color}
\usepackage{ifpdf}
\usepackage{fancybox}
\usepackage[font=small,labelfont=bf,labelsep=none]{caption}
\usepackage{float}

\allowdisplaybreaks[4]

\usepackage[latin1]{inputenc}

\newtheorem{thm}{Theorem}[section]

\newtheorem{defi}[thm]{Definition}
\newtheorem{lem}[thm]{Lemma}
\newtheorem{core}[thm]{Corollary}

\def\pf{\noindent{\it Proof.} }
\def\qed{\nopagebreak\hfill{\rule{4pt}{7pt}}
\medbreak}

\numberwithin{equation}{section}

\def\qed{\nopagebreak\hfill{\rule{4pt}{7pt}}
\medbreak}

\newlength{\boxedparwidth}
  {\begin{center} \begin{tabular}{|@{\hspace{.315in}}c@{\hspace{.15in}}|}
                  \hline \\ \begin{minipage}[t]{\boxedparwidth}
                  \setlength{\parindent}{.25in}}%
  {\end{minipage} \\ \\ \hline \end{tabular} \end{center}}

\begin{document}

\begin{center}

 {\Large \bf Some {s}eparable integer partition classes}
\end{center}

\begin{center}
 {Y.H. Chen}$^{1}$, {Thomas Y. He}$^{2}$, {F. Tang}$^{3}$ and
  {J.J. Wei}$^{4}$ \vskip 2mm

$^{1,2,3,4}$ School of Mathematical Sciences, Sichuan Normal University, Chengdu 610066, P.R. China

   \vskip 2mm

  $^1$chenyh@stu.sicnu.edu.cn, $^2$heyao@sicnu.edu.cn,  $^3$tangfan@stu.sicnu.edu.cn,  $^4$wei@stu.sicnu.edu.cn
\end{center}

{\noindent \bf Abstract.} Recently, Andrews introduced separable integer partition classes and analyzed some well-known theorems. In this paper, we investigate partitions with parts separated by parity introduced by Andrews with the aid of {separable} integer partition classes with modulus $2$. We also extend separable integer partition classes with modulus $1$ to overpartitions, called separable overpartition classes. We study overpartitions and the overpartition analogue of Rogers-Ramanujan identities, which are separable overpartition classes.

\noindent {\bf Keywords}: separable integer partition classes, overpartitions, parity of parts,  Rogers-Ramanujan identities

\noindent {\bf AMS Classifications}: 05A17, 11P84

\section{Introduction}

A partition $\pi$ of a positive integer $n$ is a finite non-increasing sequence of positive integers $\pi=(\pi_1,\pi_2,\ldots,\pi_\ell)$ such that $\pi_1+\pi_2+\cdots+\pi_\ell=n$. The $\pi_i$ are called the parts of $\pi$.  { Let $\ell(\pi)$ be the number of parts of $\pi$.} The weight of $\pi$ is the sum of parts, denoted $|\pi|$. We sometimes write  $\pi=(1^{f_1}2^{f_2}3^{f_3}\cdots)$, where $f_t(\pi)$ (or $f_t$ for short) denotes the number of parts equal to $t$ in $\pi$. { The empty sequence forms the only partition of zero.}

 In \cite{Andrews-2022}, Andrews introduced separable integer partition classes and  analyzed some well-known theorems from the point of view of separable integer partition classes, such as the first G\"ollnitz-Gordon identity, Schur's partition theorem, partitions  with $n$ copies of $n$, and so on.

 \begin{defi}
A separable integer partition class $\mathcal{P}$ with modulus $k$ is a subset of all the partitions satisfying the following{\rm:}

There is a subset $\mathcal{B}\subset\mathcal{P}$ {\rm(}$\mathcal{B}$ is called the basis of $\mathcal{P}${\rm)} such that for each integer $m\geq 1$, the number of partitions in $\mathcal{B}$ with $m$ parts is finite and every partition in $\mathcal{P}$ with $m$ parts is uniquely of the form
\begin{equation}\label{ordinary-form-1}
(b_1+\pi_1)+(b_2+\pi_2)+\cdots+(b_m+\pi_m),
\end{equation}
where $(b_1,b_2,\ldots,b_m)$ is a partition in $\mathcal{B}$ and $(\pi_1,\pi_2,\ldots,\pi_m)$ is a  non-increasing sequence of nonnegative integers, whose only restriction is that each part is divisible by $k$. Furthermore, all partitions of the form \eqref{ordinary-form-1} are in $\mathcal{P}$.
\end{defi}
For $m\geq 1$, let $b_\mathcal{B}(m)$ be the generating function for the partitions in $\mathcal{B}$ with  $m$ parts. Clearly,
\[\sum_{\pi\in\mathcal{P}}{x^{\ell(\pi)}}q^{|\pi|}=1+\sum_{m\geq 1}\frac{{x^{m}}b_\mathcal{B}(m)}{(q^k;q^k)_m}.\]

 By virtue of separable integer partition classes, Passary \cite[Section 3]{Passary-2019} studied partitions with parts separated by parity, the first little G\"ollnitz identity, the second little G\"ollnitz identity and the second G\"ollnitz-Gordon identity.

 In this paper, we investigate partitions with parts separated by parity introduced by Andrews \cite{Andrews-2018,Andrews-2019} with the aid of {separable} integer partition classes. We extend separable integer partition classes with modulus $1$ to overpartitions, called {separable} overpartition classes. An overpartition, introduced by Corteel and Lovejoy \cite{Corteel-Lovejoy-2004},  is a partition such that the first occurrence of a part can be overlined. Then, we consider the following three identities with the aid of separable overpartition classes.

\begin{equation}\label{main-1}
\sum_{n\geq0} \frac{(-a;q)_{n}q^{n}}{(q;q)_{n}}=\frac{{(-aq;q)}_{\infty}}{(q;q)_{\infty}},
\end{equation}

\begin{equation}\label{main-2}
\sum_{n\geq0} \frac{(-1;q)_{n}q^{{n+1}\choose 2}}{(q;q)_{n}}=\frac{{(-q;q^2)}_{\infty}}{(q;q^2)_{\infty}}
\end{equation}
and
\begin{equation}\label{main-3}
\sum_{n\geq0} \frac{(-q;q)_{n}q^{{n+1}\choose 2}}{(q;q)_{n}}=\frac{{(-q;q)}_{\infty}}{(q^2;q^4)_{\infty}}.
\end{equation}
Here and in the sequel, we assume that $|q|<1$ and employ the standard notation:
\[(a;q)_\infty=\prod_{i=0}^{\infty}(1-aq^i) \quad\text{and}\quad (a;q)_n=\frac{(a;q)_\infty}{(aq^n;q)_\infty}.\]

Equation \eqref{main-1} was given by Corteel and Lovejoy \cite[(1.1)]{Corteel-Lovejoy-2004}. Equations \eqref{main-2} and \eqref{main-3} were found by  \cite{Lovejoy-2003}, which can be seen as  overpartition analogues of the celebrated Rogers-Ramanujan identities, see also Chen, Sang and Shi \cite{Chen-Sang-Shi-2013}.

This paper is organized as follows. In Section 2, we recall some necessary identities. In Section 3, we investigate partitions with parts separated by parity with the aid of {separable} integer partition classes with modulus $2$. We extend separable integer partition classes with modulus $1$ to overpartitions and show \eqref{main-1} in Section 4.   Section 5 is devoted to analyzing the overpartition analogues of Rogers-Ramanujan identities, which led to proofs of \eqref{main-2} and \eqref{main-3}.

\section{Preliminaries}

In this section, we collect some well-known identities needed in this paper from \cite{Andrews-1976}.

{\noindent \bf The $q$-binomial theorem \cite[Theorem 2.1]{Andrews-1976}:}
\begin{equation}\label{chang-1}
\sum_{n\geq0} \frac{(a;q)_{n}}{(q;q)_{n}}t^{n}=\frac{{(at;q)}_{\infty}}{(t;q)_{\infty}}.
\end{equation}

It is worth mentioning that equation \eqref{main-1} can be obtained by letting $a\rightarrow -a$ and $t=q$ in \eqref{chang-1}. There are two special cases of \eqref{chang-1} (see \cite[Corollary 2.2]{Andrews-1976}):
\begin{equation}\label{Euler-1}
\sum_{n\geq0}\frac{t^n}{(q;q)_n}=\frac{1}{(t;q)_\infty},
\end{equation}
and
\begin{equation}\label{Euler-2}
\sum_{n\geq0}\frac{t^nq^{{n}\choose 2}}{(q;q)_n}=(-t;q)_\infty.
\end{equation}

{\noindent \bf Lebesgue's identity \cite[Corollary 2.7]{Andrews-1976}:}

\begin{equation}\label{chang-2}
\sum_{n\geq0} \frac{(a;q)_{n}q^{{n+1}\choose 2}}{(q;q)_{n}}={(aq;q^2)}_{\infty}(-q;q)_{\infty}.
\end{equation}

It is worth mentioning that equations \eqref{main-2} and \eqref{main-3} are the cases $a=-1$ and $a=-q$ of \eqref{chang-2} respectively.

\begin{defi}
The $q$-binomial coefficient, or Gaussian polynomial for non-negative integers $M$ and $N$ is
\[{M\brack N}=\left\{\begin{array}{ll}\frac{(q;q)_M}{(q;q)_N(q;q)_{M-N}},&\text{if }0\leq N\leq M,\\
0,&\text{otherwise.}
\end{array}
\right.\]
\end{defi}

{\noindent \bf The standard recurrences for the $q$-binomial coefficients \cite[Theorem 3.2]{Andrews-1976}:}
\begin{equation}\label{bin-r-1}
{M\brack N}={{M-1}\brack{N-1}}+q^{N}{{M-1}\brack{N}},
\end{equation}
and
\begin{equation}\label{bin-r-2}
{M\brack N}={{M-1}\brack{N}}+q^{M-N}{{M-1}\brack{N-1}}.
\end{equation}

We need the following formula related to the $q$-binomial coefficients \cite[Theorem 3.3]{Andrews-1976}.
\begin{equation}\label{j-sum}
(z;q)_{N}=\sum_{j\geq 0}{N\brack j}(-z)^{j}q^{\binom{j}{2}}.
\end{equation}

\section{Partitions with parts separated by parity}

In \cite{Andrews-2018,Andrews-2019}, Andrews considered partitions in which parts of a given parity are all smaller than those of the other parity, and if the smaller parity is odd then odd parts must appear. Andrews designated both cases, where the even (resp. odd) parts are distinct with the couplet ``$ed$" (``$od$"), or when the even (resp. odd) parts may appear an unlimited number of times with the couplet ``$eu$" (resp. ``$ou$"). The eight partition functions $p^{zw}_{xy}$  designate the partition functions in question, where $xy$ constrains the smaller parts and $zw$ the larger parts. By virtue of {separable} integer partition classes with modulus $2$, Passary studies $p^{ou}_{eu}$ and $p^{eu}_{ou}$ in \cite[Subsection 3.2]{Passary-2019}. In this section, we investigate the remaining six cases. For convenience, we use $\mathcal{P}^{zw}_{xy}$ to denote the set of partitions of type $p^{zw}_{xy}$. We will show that $\mathcal{P}^{zw}_{xy}$ is a separable integer partition class with modulus $2$.

\subsection{$p^{od}_{ed}$}

 For $m\geq1$, { let $\pi=(\pi_1,\pi_2,\ldots,\pi_m)$ be a partition in $\mathcal{P}^{od}_{ed}$. Set $j$ to be the number of even parts in $\pi$. Then, we see that $0\leq j\leq m$, $\pi_{m-j+1},\ldots,\pi_{m}$ are distinct even parts and $\pi_1,\ldots,\pi_{m-j}$ are distinct odd parts. Moreover, we have
 \[\pi_m\geq 2,\ldots,\pi_{m-j+1}\geq 2j,\]
and
 \[\pi_{m-j}\geq 2j+1,\ldots,\pi_{1}\geq 2m-1.\]
So, we can check that
\[(\pi_1-(2m-1),\ldots,\pi_{m-j}-(2j+1),\pi_{m-j+1}-2j,\ldots,\pi_m-2)\]
is a non-increasing sequence of nonnegative even integers, and vice versa.

 }

 Define
  \[\mathcal{B}^{od}_{ed}(m)=\{(2m-1,\ldots,2j+1,2j,\ldots,2)|0\leq j\leq m\}.\]

  Set
  \[\mathcal{B}^{od}_{ed}=\bigcup_{m\geq 1}\mathcal{B}^{od}_{ed}(m).\]

  Then, we see that $\mathcal{P}^{od}_{ed}$ is a separable integer partition class with modulus $2$ and $\mathcal{B}^{od}_{ed}$ is the basis of $\mathcal{P}^{od}_{ed}$. So, we get
\begin{align*}
\sum_{\pi\in \mathcal{P}^{od}_{ed}}{x^{\ell(\pi)}}q^{|\pi|}&=1+\sum_{m\geq 1}\frac{{x^m}}{(q^2;q^2)_m}\sum_{j=0}^mq^{m^2+j}\\
&=\sum_{m\geq 0}\frac{{x^m}q^{m^2}}{(q^2;q^2)_m}\frac{1-q^{m+1}}{1-q}\\
&=\frac{1}{1-q}\left(\sum_{m\geq 0}\frac{{x^m}q^{m^2}}{(q^2;q^2)_m}-q\sum_{m\geq 0}\frac{{x^m}q^{m^2+m}}{(q^2;q^2)_m}\right)\\
&=\frac{1}{1-q}\left((-{x}q;q^2)_\infty-q{(-{x}q^2;q^2)_\infty}\right),
\end{align*}
where the final equation follows from \eqref{Euler-2}. Multiplying the extremes of the foregoing string of equations by $(1-q)$, we have
\[(1-q)\sum_{\pi\in \mathcal{P}^{od}_{ed}}{x^{\ell(\pi)}}q^{|\pi|}=(-{x}q;q^2)_\infty-q{(-{x}q^2;q^2)_\infty}.\]

Let ${p}^{od}_{ed}({m,}n)$ denote the number of partitions { of $n$} in $\mathcal{P}^{od}_{ed}$ { with $m$ exactly parts} and let $D_o({m,}n)$ (resp. $D_e({m,}n)$) be the number of partitions { of $n$} with {$m$} distinct odd parts (resp. even parts). Comparing coefficients of ${x^m}q^n$ in both sides of the above identity, we deduce that

\begin{thm}
For ${m,}n\geq1,$ we have
\[{p}^{od}_{ed}({m,}n)-{p}^{od}_{ed}({m,}n-1)=D_o({m,}n)-D_e({m,}n-1).\]
\end{thm}

It is easy to see that $D_e({m,}n-1)=0$ when $n$ is even. Thus, we derive that
\begin{core}
For ${m,}n\geq1,$ we have
\[{p}^{od}_{ed}({m,}2n)-{p}^{od}_{ed}({m,}2n-1)=D_o({m,}2n).\]
\end{core}

{ Moreover, since $D_o(m,2n)=0$ when $m$ is odd, we get
\begin{core}
For $m,n\geq1,$ we have
\[{p}^{od}_{ed}(2m-1,2n)={p}^{od}_{ed}(2m-1,2n-1).\]
\end{core}}

\subsection{$p^{ed}_{od}$}

 For $m\geq1$, { let $\pi=(\pi_1,\pi_2,\ldots,\pi_m)$ be a partition in $\mathcal{P}^{ed}_{od}$. Set $j$ to be the number of odd parts in $\pi$. Recall that  odd parts must appear if the smaller parity is odd. Then, odd parts must appear in $\pi$, and so we have $1\leq j\leq m$. We see that $\pi_{m-j+1},\ldots,\pi_{m}$ are distinct odd parts  and $\pi_1,\ldots,\pi_{m-j}$ are distinct even parts. Moreover, we have
 \[\pi_m\geq 1,\ldots,\pi_{m-j+1}\geq 2j-1,\]
and
 \[\pi_{m-j}\geq 2j,\ldots,\pi_{1}\geq 2m-2.\]
So, we can check that
\[(\pi_1-(2m-2),\ldots,\pi_{m-j}-2j,\pi_{m-j+1}-(2j-1),\ldots,\pi_m-1)\]
is a non-increasing sequence of nonnegative even integers, and vice versa.

 }

 Define
 \[\mathcal{B}^{ed}_{od}(m)=\{(2m-2,\ldots,2j,2j-1,\ldots,1)|1\leq j\leq m\}.\]
  Set
  \[\mathcal{B}^{ed}_{od}=\bigcup_{m\geq 1}\mathcal{B}^{ed}_{od}(m).\]

Then, we see that $\mathcal{P}^{ed}_{od}$ is a separable integer partition class with modulus $2$ and $\mathcal{B}^{ed}_{od}$ is the basis of $\mathcal{P}^{ed}_{od}$. So, we get
\begin{align*}
\sum_{\pi\in \mathcal{P}^{ed}_{od}}{x^{\ell(\pi)}}q^{|\pi|}&=\sum_{m\geq 1}\frac{{x^m}}{(q^2;q^2)_m}\sum_{j=1}^mq^{m^2-m+j}\\
&=\sum_{m\geq 0}\frac{{x^m}q^{m^2-m+1}}{(q^2;q^2)_m}\frac{1-q^{m}}{1-q}\\
&=\frac{q}{1-q}\left(\sum_{m\geq 0}\frac{{x^m}q^{m^2-m}}{(q^2;q^2)_m}-\sum_{m\geq 0}\frac{{x^m}q^{m^2}}{(q^2;q^2)_m}\right)\\
&=\frac{q}{1-q}\left({(-x;q^2)_\infty}-{(-{x}q;q^2)_\infty}\right),
\end{align*}
where the final equation follows from \eqref{Euler-2}. Multiplying the extremes of the foregoing string of equations by $(1-q)$, we have
\[(1-q)\sum_{\pi\in \mathcal{P}^{ed}_{od}}{x^{\ell(\pi)}}q^{|\pi|}=q\left({(-x;q^2)_\infty}-{(-{x}q;q^2)_\infty}\right).\]

Let ${p}^{ed}_{od}({m,}n)$ denote the number of partitions { of $n$} in $\mathcal{P}^{ed}_{od}$ { with exactly $m$ parts}. Comparing coefficients of ${x^m}q^n$ in  both sides of the above identity, we deduce that

\begin{thm}\label{thm-od-ed}
For ${m,}n\geq1,$ we have
\begin{equation*}\label{sum-with-m}
{p}^{ed}_{od}({m,}n)-{p}^{ed}_{od}({m,}n-1)={D_e(m,n-1)+D_e(m-1,n-1)-D_o(m,n-1).}
\end{equation*}
\end{thm}

It yields that
{\begin{core}
For $n\geq1,$ we have
\[{p}^{ed}_{od}(n)-{p}^{ed}_{od}(n-1)=2D_e(n-1)-D_o(n-1),\]
and
\[{p}^{ed}_{od}(n)-{p}^{ed}_{od}(n-1)\equiv D_o(n-1)\pmod{2}.\]
\end{core}
}

Again by $D_e({m,}n-1)=0$ when $n$ is even, it follows from Theorem \ref{thm-od-ed} that
\begin{core}
For ${m,} n\geq 1,$ we have
\[{p}^{ed}_{od}({m,}2n)-{p}^{ed}_{od}({m,}2n-1)=-D_o({m,}2n-1).\]
\end{core}

{ Moreover, since $D_o(m,2n-1)=0$ when $m$ is even, we get
\begin{core}
For $m,n\geq1,$ we have
\[{p}^{ed}_{od}(2m,2n)={p}^{ed}_{od}(2m,2n-1).\]
\end{core}}

\subsection{The remaining four cases}

Of the remaining four cases, we can show that they are also separable integer partition classes with modulus $2$ {by using the similar argument in Subsections 3.1 and 3.2. We will omit the proofs here.}  For $m\geq 1$, we use $\mathcal{B}^{zw}_{xy}(m)$ to denote the set of partitions in the basis of  $\mathcal{P}^{zw}_{xy}$ with exactly $m$ parts. We just present the $\mathcal{B}^{zw}_{xy}(m)$ and the generating function for the partitions in $\mathcal{P}^{zw}_{xy}$. {But the results can not be simplified.} In this subsection, we write partition of the form $(1^{f_1}2^{f_2}3^{f_3}\cdots)$.

{\noindent \bf $\mathcal{P}^{ed}_{ou}$:}  For $m\geq1$, define
 \[\mathcal{B}^{ed}_{ou}(m)=\{(1^{j},2,4,\ldots,2m-2j)|1\leq j\leq m\}.\]

The generating function for the partitions in $\mathcal{P}^{ed}_{ou}$ is
\begin{align*}
\sum_{\pi\in \mathcal{P}^{ed}_{ou}}{x^{\ell(\pi)}}q^{|\pi|}&=\sum_{m\geq 1}\frac{{x^m}}{(q^2;q^2)_m}\sum_{\pi\in \mathcal{B}^{ed}_{ou}(m)}q^{|\pi|}\\
&=\sum_{m\geq 1}\sum_{j=1}^m\frac{{x^m}q^{m^2-2mj+j^2+m}}{(q^2;q^2)_m}.
\end{align*}

{\noindent \bf $\mathcal{P}^{eu}_{od}$:}  For $m\geq1$, define
 \[\mathcal{B}^{eu}_{od}(m)=\{(1,3,\ldots,2j-1,2j^{m-j})|1\leq j\leq m\}.\]

The generating function for the partitions in $\mathcal{P}^{eu}_{od}$ is
\begin{align*}
\sum_{\pi\in \mathcal{P}^{eu}_{od}}{x^{\ell(\pi)}}q^{|\pi|}&=\sum_{m\geq 1}\frac{{x^m}}{(q^2;q^2)_m}\sum_{\pi\in \mathcal{B}^{eu}_{od}(m)}q^{|\pi|}\\
&=\sum_{m\geq 1}\sum_{j=1}^m\frac{{x^m}q^{2mj-j^2}}{(q^2;q^2)_m}.
\end{align*}

{\noindent \bf $\mathcal{P}^{od}_{eu}$:}  For $m\geq1$, define
 \[\mathcal{B}^{od}_{eu}(m)=\{(2^{j},3,5,\ldots,2m-2j+1)|1\leq j\leq m\}\bigcup{\{(1,3,\ldots,2m-1)\}}.\]

The generating function for the partitions in $\mathcal{P}^{ed}_{ou}$ is
\begin{align*}
\sum_{\pi\in \mathcal{P}^{od}_{eu}}{x^{\ell(\pi)}}q^{|\pi|}&=1+\sum_{m\geq 1}\frac{{x^m}}{(q^2;q^2)_m}\sum_{\pi\in \mathcal{B}^{od}_{eu}(m)}q^{|\pi|}\\
&=\sum_{m\geq 1}\sum_{j=1}^m\frac{{x^m}q^{m^2-2mj+j^2+2m}}{(q^2;q^2)_m}+\sum_{m\geq 0}\frac{{x^m}q^{m^2}}{(q^2;q^2)_m}\\
&=\sum_{m\geq 1}\sum_{j=1}^m\frac{{x^m}q^{m^2-2mj+j^2+2m}}{(q^2;q^2)_m}+(-{x}q;q^2)_\infty,
\end{align*}
where the final equation follows from \eqref{Euler-2}.

{\noindent \bf $\mathcal{P}^{ou}_{ed}$:}  For $m\geq1$, define
 \[{\mathcal{B}^{ou}_{ed}(m)}=\{(2,4,\ldots,2j,(2j+1)^{m-j})|1\leq j\leq m\}\bigcup{\{(1^m)\}}.\]

The generating function for the partitions in $\mathcal{P}^{ou}_{ed}$ is
\begin{align*}
\sum_{\pi\in \mathcal{P}^{ou}_{ed}}{x^{\ell(\pi)}}q^{|\pi|}&=1+\sum_{m\geq 1}\frac{{x^m}}{(q^2;q^2)_m}\sum_{\pi\in {\mathcal{B}^{ou}_{ed}(m)}}q^{|\pi|}\\
&=\sum_{m\geq 1}\sum_{j=1}^m\frac{{x^m}q^{2mj-j^2+m}}{(q^2;q^2)_m}+\sum_{m\geq 0}\frac{{x^m}q^{m}}{(q^2;q^2)_m}\\
&=\sum_{m\geq 1}\sum_{j=1}^m\frac{{x^m}q^{2mj-j^2+m}}{(q^2;q^2)_m}+\frac{1}{({x}q;q^2)_\infty},
\end{align*}
where the final equation follows from \eqref{Euler-1}.

\section{Overparitions}

Let $\overline{\mathcal{P}}$ denote the set of all overpartitions. In this section, we will investigate the set $\overline{\mathcal{P}}$ and then give a proof of \eqref{main-1}. We impose the following order on the parts of an overpartition:
\[1<\bar{1}<2<\bar{2}<\cdots.\]
We adopt the following convention: For positive integers $t$ and $b$, we say that a part is of size $t$ if the part is $t$ or $\bar{t}$ and we define $t\pm b$ (resp. $\overline{t}\pm b$) as a non-overlined part (resp. an overlined part) of size  $t\pm b$.

\subsection{Separable overpartition classes}

We extend separable integer partition classes with modulus $1$ to overpartitions, called separable overpartition classes, which is stated as follows.

\begin{defi}
A separable overpartition class $\mathcal{P}$ is a subset of all the overpartitions satisfying the following{\rm:}

There is a subset $\mathcal{B}\subset\mathcal{P}$ {\rm(}$\mathcal{B}$ is called the basis of $\mathcal{P}${\rm)} such that for each integer $m\geq 1$, the number of overpartitions in $\mathcal{B}$ with $m$ parts is finite and every overpartition in $\mathcal{P}$ with $m$ parts is uniquely of the form
\begin{equation}\label{over-form-1}
(b_1+\pi_1)+(b_2+\pi_2)+\cdots+(b_m+\pi_m),
\end{equation}
where $(b_1,b_2,\ldots,b_m)$ is an overpartition in $\mathcal{B}$ and $(\pi_1,\pi_2,\ldots,\pi_m)$ is a  non-increasing sequence of nonnegative integers. Moreover, all overpartitions of the form \eqref{over-form-1} are in $\mathcal{P}$.
\end{defi}

Let $b_\mathcal{B}(m)$ be the generating function for the overpartitions in $\mathcal{B}$ with  $m$ parts and let $b_\mathcal{B}(m,j)$ (resp. $b_\mathcal{B}(m,\bar{j})$) be the generating function for the overpartitions in $\mathcal{B}$ with exactly $m$ parts and { the} largest part $j$ (resp. $\bar{j}$). Clearly,
\[b_\mathcal{B}(m)=\sum_{j\geq1}\left[b_\mathcal{B}(m,j)+b_\mathcal{B}(m,\bar{j})\right].\]
 The generating function for  the overpartitions in $\mathcal{P}$ is given by
 \[\begin{split}\sum_{\pi\in\mathcal{P}}q^{|\pi|}&=\sum_{m\geq0}\frac{b_\mathcal{B}(m)}{(q;q)_m}\\
 &=1+\sum_{m\geq1}\frac{1}{(q;q)_m}\sum_{j\geq1}\left[b_\mathcal{B}(m,j)+b_\mathcal{B}(m,\bar{j})\right].
 \end{split}\]

\subsection{The basis of $\overline{\mathcal{P}}$}

 We aim to show that $\overline{\mathcal{P}}$ is a separable overpartition class. To do this, we are required to find the basis of $\overline{\mathcal{P}}$, which involves the following set.
\begin{defi}\label{base-over}
For $m\geq 1$, define $\mathcal{O}(m)$ to be the set of overpartitions $\pi=(\pi_{1},\pi_{2},\cdots \pi_{m})$ such that
\begin{itemize}
\item[{\rm(1)}] $\pi_{m}=1$ or $\bar{1}${\rm;}
\item[{\rm(2)}] for $1\leq i<m,$ $\pi_{i}\leq\pi_{i+1}+1$ with strict inequality if $\pi_{i+1}$ is non-overlined.
\end{itemize}
\end{defi}

For $m\geq 1$, assume that $\pi$ is an overpartition in {$\mathcal{O}(m)$}. For $1\leq i<m$, if $\pi_{i+1}=t$, then we have $t\leq\pi_{i}<t+1$, and so  $\pi_{i}=t$ or $\bar{t}$. If $\pi_{i+1}=\bar{t}$, then we have $\overline{t}<\pi_{i}\leq \overline{t+1}$, and so  $\pi_{i}=t+1$ or $\overline{t+1}$.
{Then, we see that there are $2^m$ overpartitions in $\mathcal{O}(m)$.

For example, the number of overpartitions in $\mathcal{O}(4)$ is $2^4=16$.
\[(1,1,1,1),(\bar{1},1,1,1),(2,\bar{1},1,1),(\bar{2},\bar{1},1,1),(2,2,\bar{1},1),(\bar{2},2,\bar{1},1),(3,\bar{2},\bar{1},1),(\bar{3},\bar{2},\bar{1},1),\]
\[(2,2,2,\bar{1}),(\bar{2},2,2,\bar{1}),(3,\bar{2},2,\bar{1}),(\bar{3},\bar{2},2,\bar{1}),(3,3,\bar{2},\bar{1}),(\bar{3},3,\bar{2},\bar{1}),(4,\bar{3},\bar{2},\bar{1}),(\bar{4},\bar{3},\bar{2},\bar{1}).\]

}

\begin{lem}\label{over-lem-1}
 $\overline{\mathcal{P}}$ is a separable overpartition class.
\end{lem}

\pf Set
\[\mathcal{O}=\bigcup_{m\geq1}\mathcal{O}(m).\]
 Obviously, $\mathcal{O}$ is the basis of $\overline{\mathcal{P}}$. This completes the proof.   \qed

Let $b_\mathcal{O}(m,j)$ (resp. $b_\mathcal{O}(m,\bar{j})$) be the generating function for the overpartitions in $\mathcal{O}(m)$ with the largest part $j$ (resp. $\bar{j}$). Then, we have
\begin{lem}\label{over-lem-2}
For $m\geq 1$ and $j\geq 1$, we have
\[
b_\mathcal{O}(m,\bar{j})=b_\mathcal{O}(m,j)=q^{m+\binom{j}{2}}{{m-1}\brack{j-1}}.
\]
\end{lem}
\pf It suffices to show that
\begin{itemize}
\item[(A)] $b_\mathcal{O}(m,\bar{j})=b_\mathcal{O}(m,j)$;
\item[(B)] $b_\mathcal{O}(m,j)=q^{m+\binom{j}{2}}{{m-1}\brack{j-1}}$.
\end{itemize}

{\noindent Condition (A).} For an overpartition $\pi$ in $\mathcal{O}(m)$ with the largest part $\bar{j}$, we can get an overpartition in $\mathcal{O}(m)$ with the largest part ${j}$ by changing the overlined part $\bar{j}$ in $\pi$ to a non-overlined part $j$, and vice {versa}.
This implies that Condition (A) is satisfied.

{\noindent Condition (B).} We proceed to show Condition (B) by induction on $m$. For $m=1$, it is easy to check that
\begin{eqnarray*}
b_\mathcal{O}(1,j)=
\begin{cases}
q, \quad \text{if} \ j=1, \\[5pt]
0, \quad \text{otherwise.}
\end{cases}
\end{eqnarray*}
So, Condition (B) holds for $m=1$.

Assume that  $m\geq 2$ and Condition (B) holds for $m-1$, that is, for $j\geq 1$,
\begin{equation}\label{ying-1}
b_\mathcal{O}(m-1,j)=q^{m-1+\binom{j}{2}}{{m-2}\brack{j-1}}.
\end{equation}

We proceed to show that Condition (B) also holds for $m$. The only one overpartition in $\mathcal{O}(m)$ with the largest part $1$ is $(1^m)$.
This implies that   $b_\mathcal{O}(m,1)=q^m$, which agrees with Condition (B) for $j=1$. For $j\geq 2$, we need to construct the following relation.
\begin{equation}\label{relation-1}
b_\mathcal{O}(m,j)=q^{j}\left[B(m-1,j)+B(m-1,j-1)\right].
\end{equation}

Assume that $\pi=(\pi_1,\pi_2,\ldots,\pi_m)$ is an {overpartition in} $\mathcal{O}(m)$ with the largest part ${j}$. Then, we have $\pi_1=j$. By the condition (2) in Definition \ref{base-over}, we deduce that $\pi_2=j$ or $\overline{j-1}$. If we remove the largest part $j$ from $\pi$, then we can get an overpartition in $\mathcal{O}(m-1)$ with the largest part $j$ or $\overline{j-1}$. It yields
\begin{equation*}\label{relation-2}
b_\mathcal{O}(m,j)=q^{j}\left[B(m-1,j)+B(m-1,\overline{j-1})\right].
\end{equation*}
Combining with Condition (A), we arrive at \eqref{relation-1}. Using \eqref{ying-1} and \eqref{relation-1}, we get
\begin{align*}
b_\mathcal{O}(m,j)&=q^{j}\left(q^{m-1+\binom{j}{2}}{{m-2}\brack{j-1}}+
q^{m-1+\binom{j-1}{2}}{{m-2}\brack{j-2}}\right)\\[5pt]
&=q^{m+\binom{j}{2}}\left(q^{j-1}{{m-2}\brack{j-1}}+
{{m-2}\brack{j-2}}\right)\\[5pt]
&=q^{m+\binom{j}{2}}{{m-1}\brack{j-1}},
\end{align*}
where the final equation follows from \eqref{bin-r-1}. So, Condition (B) is valid for $j\geq 2$. We conclude that Condition (B) also holds for $m$. This completes the proof.  \qed

\subsection{Proof of \eqref{main-1}}

For an overpartition $\pi$, let $\ell_o(\pi)$ be the number of overlined parts in $\pi$. In order to give a proof of \eqref{main-1}, we need the following Lemma.
\begin{lem}
\begin{equation}\label{over-main-identity}
\sum_{\pi\in{\overline{\mathcal{P}}}}a^{\ell_o(\pi)}q^{|\pi|}=\sum_{m\geq0}\frac{q^m}{(q;q)_m}\sum_{j\geq 0}a^jq^{j\choose 2}{m\brack j}.
\end{equation}
\end{lem}

\pf Assume that  $\pi$ is an overpartition in $\mathcal{O}$. If the largest part of $\pi$ is $\overline{j}$, then we have $\ell_o(\pi)=j$ since $\overline{1}$, $\overline{2},\ldots$, $\overline{j}$ are parts of $\pi$. If the largest part of $\pi$ is ${j}$, then we have $\ell_o(\pi)=j-1$ because $\overline{1}$, $\overline{2},\ldots$, $\overline{j-1}$ are parts of $\pi$. Then, appealing to Lemmas \ref{over-lem-1} and \ref{over-lem-2}, we get
\begin{align*}
\sum_{\pi\in{\overline{\mathcal{P}}}}a^{\ell_o(\pi)}q^{|\pi|}&=1+\sum_{m\geq1}\frac{1}{(q;q)_m}\sum_{j\geq0}\left[a^{j}b_{\mathcal{O}}(m,j+1)+a^{j}b_{\mathcal{O}}(m,\bar{j})\right]\\[5pt]
&=1+\sum_{m\geq 1}\frac {1}{(q;q)_{m}}\sum_{j\geq 0}\left(a^jq^{m+\binom{j+1}{2}}
{{m-1}\brack j}+a^jq^{m+\binom{j}{2}}{{m-1}\brack {j-1}}\right)\\[5pt]
&=1+\sum_{m\geq 1}\frac {q^m}{(q;q)_{m}}\sum_{j\geq 0}a^jq^{\binom{j}{2}}\left(q^{j}
{{m-1}\brack j}+{{m-1}\brack {j-1}}\right)\\[5pt]
&=1+\sum_{m\geq 1}\frac {q^m}{(q;q)_{m}}\sum_{j\geq 0}a^jq^{\binom{j}{2}}
{{m}\brack j},
\end{align*}
where the final equation follows from \eqref{bin-r-1}. This completes the proof.  \qed

  We are now in a position to give a proof of \eqref{main-1}.

{\noindent \bf Proof of \eqref{main-1}.} We first sum up the $j$-sum in the right-hand side of \eqref{over-main-identity} by letting $z=-a$ and $N=m$ in \eqref{j-sum}. So, we get
\begin{equation}\label{over-main-identity-1}
\sum_{\pi\in{\overline{\mathcal{P}}}}a^{\ell_o(\pi)}q^{|\pi|}=\sum_{m\geq0}\frac{(-a;q)_mq^m}{(q;q)_m}.
\end{equation}

On the other hand, by interchanging the order of summation in the right-hand side of \eqref{over-main-identity}, we can obtain that
\begin{equation}\label{over-main-identity-2}
\begin{split}
\sum_{\pi\in{\overline{\mathcal{P}}}}a^{\ell_o(\pi)}q^{|\pi|}&=\sum_{j\geq0}\frac {a^jq^{\binom{j}{2}}}{(q;q)_{j}}\sum_{m\geq j}\frac {q^{m}}{(q;q)_{m-j}}\\
&=\sum_{j\geq0}\frac {a^jq^{\binom{j+1}{2}}}{(q;q)_{j}}\sum_{m\geq 0}\frac {q^{m}}{(q;q)_{m}}\\
&=\frac {(-aq;q)_{\infty}}{(q;q)_{\infty}},
\end{split}
\end{equation}
where the final equation follows from \eqref{Euler-1} and \eqref{Euler-2}. Then, \eqref{main-1} immediately follows from \eqref{over-main-identity-1} and \eqref{over-main-identity-2}. Thus, we complete the proof.   \qed

\section{Overpartition analogue of Rogers-Ramanujan identities}

In this section, we will investigate the overpartition analogue of Rogers-Ramanujan identities. For $r=1$ or $2$, let $\mathcal{C}_r$ denote the set of overpartitions $\pi=(\pi_1,\pi_2,\ldots,\pi_\ell)$ such that

\begin{itemize}
\item[(1)] for $1\leq i<\ell$, $\pi_{i}\geq\pi_{i+1}+1$ with strict inequality if $\pi_i$ is non-overlined;

\item[(2)] at most $r-1$ parts of the $\pi_i$ are equal to $1$.
\end{itemize}

Bear in mind that the parts in an overpartition are ordered as follows:
\begin{equation*}
1<\bar{1}<2<\bar{2}<\cdots.
\end{equation*}
 For positive integers $t$ and $b$,  we define $t\pm b$ (resp. $\overline{t}\pm b$) as a non-overlined part of size  $t\pm b$ (resp. an overlined part of size ${t\pm b}$).

Then we will study the sets $\mathcal{C}_2$ and $\mathcal{C}_1$ from the point of view of separable overpartition classes, and then give proofs of \eqref{main-2} and $\eqref{main-3}$ respectively.

 \subsection{The basis of $\mathcal{C}_r$}

For $r=1$ or $2$, we aim to show that $\mathcal{C}_r$ is a  separable overpartition class. To do this, we are obliged to give the basis of $\mathcal{C}_r$, which involves the following set.

\begin{defi}\label{base-r-r}
For $m\geq 1$, define $\mathcal{G}(m)$ to be the set of overpartitions $\pi=(\pi_{1},\pi_{2},\cdots \pi_{m})$  such that for $1\leq i<m,$
\begin{itemize}

\item[{\rm(1)}] $\pi_{i}\geq\pi_{i+1}+1$ with strict inequality if $\pi_i$ is non-overlined{\rm;}

\item[{\rm(2)}] $\pi_{i}\leq\pi_{i+1}+2$ with strict inequality if $\pi_{i}$ is overlined.
 \end{itemize}
\end{defi}

For $m\geq 1$, assume that $\pi$ is an overpartition in $\mathcal{G}(m)$. For $1\leq i<m$, if $\pi_i=t$, then we have ${t-2}\leq\pi_{i+1}<t-1$, and so  $\pi_{i+1}=t-2$ or $\overline{t-2}$. If $\pi_i=\bar{t}$, then we have $\overline{t-2}<\pi_{i+1}\leq \overline{t-1}$, and so  $\pi_{i+1}=t-1$ or $\overline{t-1}$.

{ In other words, for $m\geq 1$, assume that $\pi$ is an overpartition in $\mathcal{G}(m)$. For $1\leq i<m$, if $\pi_{i+1}$ is of size $t$, then we have $\pi_i=\overline{t+1}$ or $t+2$. Then, we see that there are $2^{m-1}$ overpartitions in $\mathcal{G}(m)$ with fixed smallest part.
For example, the number of overpartitions in $\mathcal{G}(4)$ with the smallest part $1$ or $\bar{1}$ or $2$ is $2^{4-1}\times 3=24$.
{\small\begin{equation}\label{exam-g-1}
(\bar{4},\bar{3},\bar{2},1),(5,\bar{3},\bar{2},1),(\bar{5},4,\bar{2},1),(6,4,\bar{2},1),(\bar{5},\bar{4},3,1),(6,\bar{4},3,1),(\bar{6},5,3,1),(7,5,3,1),
\end{equation}
\begin{equation}\label{exam-g-2}
(\bar{4},\bar{3},\bar{2},\bar{1}),(5,\bar{3},\bar{2},\bar{1}),(\bar{5},4,\bar{2},\bar{1}),(6,4,\bar{2},\bar{1}),(\bar{5},\bar{4},3,\bar{1}),(6,\bar{4},3,\bar{1}),(\bar{6},5,3,\bar{1}),(7,5,3,\bar{1}),
\end{equation}
\begin{equation}\label{exam-g-3}
(\bar{5},\bar{4},\bar{3},2),(6,\bar{4},\bar{3},2),(\bar{6},5,\bar{3},2),(7,5,\bar{3},2),(\bar{6},\bar{5},4,2),(7,\bar{5},4,2),(\bar{7},6,4,2),(8,6,4,2).
\end{equation}
}
}

For $m\geq 1$, define $\mathcal{G}_2(m)$ (resp. $\mathcal{G}_1(m)$) to be the set of overpartitions in $\mathcal{G}(m)$ with the smallest part $1$ or $\bar{1}$ (resp. with the smallest part  $\bar{1}$ or $2$).
{ Then, we see that there are $2^m$ overpartitions in $\mathcal{G}_2(m)$ (resp. $\mathcal{G}_1(m)$). For example, the number of overpartitions in $\mathcal{G}_2(4)$ (resp. $\mathcal{G}_1(4)$) is $2^4=16$, and
the overpartitions in $\mathcal{G}_2(4)$ (resp. $\mathcal{G}_1(4)$) are given in \eqref{exam-g-1} and \eqref{exam-g-2} (resp.  \eqref{exam-g-2} and \eqref{exam-g-3}).}

For $r=1$ or $2$, set
\[\mathcal{G}_r=\bigcup_{m\geq 1}\mathcal{G}_r(m),\]
it is easy to check that  $\mathcal{G}_r$ is the basis of $\mathcal{C}_r$. So, we have

\begin{lem}\label{R-R-lem-1}
For $r=1$ or $2,$ ${\mathcal{C}}_r$ is a separable overpartition class.
\end{lem}

For $r=1$ or $2$, the generating function for the overpartitions in $\mathcal{C}_r$ is
\[\sum_{\pi\in\mathcal{C}_r}q^{|\pi|}=1+\sum_{m\geq 1}\frac{1}{(q;q)_m}\sum_{\pi\in\mathcal{G}_r(m)}q^{|\pi|}{.}\]

Then, we proceed to investigate the generating function for the overpartitions in $\mathcal{G}_2(m)$ and $\mathcal{G}_1(m)$, and then give proofs of \eqref{main-2}
and \eqref{main-3} respectively.

\subsection{Overpartition analogue of the first Rogers-Ramanujan identity}

In this subsection, we focus on $\mathcal{G}_2(m)$ and  the proof of \eqref{main-2}. Let $b_\mathcal{G}(m,j)$ (resp. $b_\mathcal{G}(m,\bar{j})$) be the generating function for the overpartitions in $\mathcal{G}_2(m)$ with the largest part $j$ (resp. $\bar{j}$). Then, we have

\begin{lem}\label{g-1-lemmm-1}
For $m\geq 2$ and $j\geq 2,$
\begin{equation}\label{g-1-lem-1}
b_\mathcal{G}(m,j)=qb_\mathcal{G}(m,\overline{j-1}).
\end{equation}
\end{lem}

\pf Assume that $\pi$ is an overpartition in $\mathcal{G}_2(m)$ with the largest part $j$. If we change the largest part $j$ in $\pi$ to an overlined part $\overline{j-1}$, then we obtain an overpartition in $\mathcal{G}_2(m)$ with the largest part $\overline{j-1}$, and vice {versa}. This implies that \eqref{g-1-lem-1} is valid. Thus, the proof is complete.   \qed

So, we just need to give $b_\mathcal{G}(m,\bar{j})$. To do this, we need to show the following recurrence for $b_\mathcal{G}(m,\bar{j})$.

\begin{lem}\label{le-ref}
For $m\geq 3$ and $j\geq 3,$
\begin{equation}\label{g-1-lem-2}
b_\mathcal{G}(m,\overline{j})=q^{j}\left[qb_\mathcal{G}(m-1,\overline{j-2})+b_\mathcal{G}(m-1,\overline{j-1})\right].
\end{equation}
\end{lem}

\pf Assume that $\pi=(\pi_1,\pi_2,\ldots,\pi_m)$ is an overpartition in $\mathcal{G}_2(m)$ with the largest part $\bar{j}$, that is, $\pi_1=\bar{j}$. Moreover, we have $\pi_2=j-1$ or $\overline{j-1}$. If we remove the part $\bar{j}$ from $\pi$, then we get an overpartition in $\mathcal{G}_2(m-1)$ with the largest part $j-1$ or $\overline{j-1}$, and vice {versa}. It implies that
\begin{equation}\label{g-1-lem-3}
b_\mathcal{G}(m,\overline{j})=q^{j}\left[b_\mathcal{G}(m-1,{j-1})+b_\mathcal{G}(m-1,\overline{j-1})\right].
\end{equation}

Using \eqref{g-1-lem-1}, we find that
\begin{equation*}
b_\mathcal{G}(m-1,j-1)=qb_\mathcal{G}(m-1,\overline{j-2}).
\end{equation*}
Combining with \eqref{g-1-lem-3}, we arrive at  \eqref{g-1-lem-2}. Thus, we complete the proof.  \qed

Now, we are in a position to give $b_\mathcal{G}(m,\bar{j})$.
\begin{lem}\label{g-1-lemmm-2}
For $m\geq 2$ and $j\geq 0$,
\begin{equation}\label{rr-1-lem-pr}
b_\mathcal{G}(m,\overline{m+j})=2q^{\binom{m+1}{2}+\frac{j(j+3)}{2}}{{m-2}\brack j}.
\end{equation}
\end{lem}

\pf We prove this lemma by induction on $m$. For $m=2$, there are two overpartitions in $\mathcal{G}_2(m)$ such that the largest part is overlined, which are $(\bar{2},1)$ and $(\bar{2},\bar{1})$. This implies that
\begin{eqnarray*}
b_\mathcal{G}(2,\overline{2+j})=\left\{
\begin{array}{ll}
2q^{3},&\text{if }j=0,\\[5pt]
0,&\text{otherwise.}
\end{array}\right.
\end{eqnarray*}

For $m\geq 3$, assume that \eqref{rr-1-lem-pr} is valid for $m-1$, that is, for $j\geq 0$,
\[b_\mathcal{G}(m-1,\overline{m+j-1})=2q^{\binom{m}{2}+\frac{j(j+3)}{2}}{{m-3}\brack j}.\]
Combining with Lemma \ref{le-ref}, we get
\begin{align*}
b_\mathcal{G}(m,\overline{m+j})&=q^{m+j}(qb_\mathcal{G}(m-1,\overline{m+j-2})+b_\mathcal{G}(m-1,\overline{m+j-1}))\\[5pt]
&=2q^{m+j}\left(q^{1+\binom{m}{2}+\frac{(j-1)(j+2)}{2}}{{m-3}\brack {j-1}}
+q^{\binom{m}{2}+\frac{j(j+3)}{2}}{{m-3}\brack {j}}\right)\\[5pt]
&=2q^{\binom{m+1}{2}+\frac{j(j+3)}{2}}\left({{m-3}\brack {j-1}}
+q^{j}{{m-3}\brack {j}}\right)\\[5pt]
&=2q^{\binom{m+1}{2}+\frac{j(j+3)}{2}}{{m-2}\brack j},
\end{align*}
where the final equation follows from \eqref{bin-r-1}. So, conclude that  \eqref{rr-1-lem-pr} holds for $m$. This completes the proof.  \qed

For $m\geq 1$ and $j\geq 1$, let $b_\mathcal{G}(m,j,1)$ (resp. $b_\mathcal{G}(m,j,\bar{1})$) be the generating function for the overpartitions in $\mathcal{G}(m)$ with the largest part of size $j$ and the smallest part $1$ (resp. the smallest part $\bar{1}$). Clearly,
\[b_\mathcal{G}(m,j)+b_\mathcal{G}(m,\bar{j})=b_\mathcal{G}(m,j,1)+b_\mathcal{G}(m,j,\bar{1}).\]

Then, we have
\begin{lem}\label{rr-lem-relation} For $m\geq 1$ and $j\geq 1$,
\begin{equation}\label{rr-relation-1}
b_\mathcal{G}(m,j,1)=b_\mathcal{G}(m,j,\bar{1}).
\end{equation}
\end{lem}
\pf Assume that $\pi$ is an overpartitions in $\mathcal{G}(m)$ with the largest part of size $j$ and the smallest part $1$. If we change the smallest part $1$ in $\pi$ to an overlined part $\bar{1}$, then we obtain an overpartition in $\mathcal{G}(m)$ with the largest part of size $j$ and the smallest part $\bar{1}$, and vice {versa}. This implies that \eqref{rr-relation-1} is valid. Thus, the proof is complete.   \qed

\begin{lem}\label{rr-1-new}
For $m\geq 1$ and $j\geq 0$,
\begin{equation}\label{rr-2-pr}
b_\mathcal{G}(m,m+j,1)=b_\mathcal{G}(m,m+j,\bar{1})=q^{\binom{m+1}{2}+\binom{j+1}{2}}
{{m-1}\brack j}.
\end{equation}
\end{lem}

\pf  For $m=1$, there are two overpartitions $(1)$ and $(\bar{1})$ in $\mathcal{G}(m)$, which implies that
\begin{eqnarray*}
b_\mathcal{G}(1,{1+j},1)=b_\mathcal{G}(1,{1+j},\bar{1})=\left\{
\begin{array}{cc}
q,&\text{if }j=0,\\[5pt]
0,&\text{otherwise.}
\end{array}\right.
\end{eqnarray*}

For $m\geq 2$, by Lemmas \ref{g-1-lemmm-1}, \ref{g-1-lemmm-2} and \ref{rr-lem-relation}, we get
\begin{align*}
&\quad b_\mathcal{G}(m,m+j,1)=b_\mathcal{G}(m,m+j,\bar{1})\\[5pt]
&=\left[b_\mathcal{G}(m,\overline{m+j})+b_\mathcal{G}(m,m+j)\right]/2\\[5pt]
&=\left[b_\mathcal{G}(m,\overline{m+j})+qb_\mathcal{G}(m,\overline{m+j-1})\right]/2\\[5pt]
&=q^{\binom{m+1}{2}}
\left(q^{\frac{j(j+3)}{2}}{{m-2}\brack j}
+q^{1+\frac{(j-1)(j+2)}{2}}{{m-2}\brack {j-1}}\right)\\[5pt]
&=q^{\binom{m+1}{2}+\binom{j+1}{2}}\left(q^j{{m-2}\brack j}
+{{m-2}\brack {j-1}}\right)\\[5pt]
&=q^{\binom{m+1}{2}+\binom{j+1}{2}}
{{m-1}\brack j},
\end{align*}
where the final equation follows from \eqref{bin-r-1}. This completes the proof.   \qed

Now, we give the generating function for the overpartitions in $\mathcal{C}_2$.
\begin{lem}
\begin{equation}\label{rr-1-main}
\sum_{\pi\in\mathcal{C}_2}q^{|\pi|}=\sum_{m\geq0}\frac{q^{\binom{m+1}{2}}}{(q;q)_{m}}
\sum_{j\geq0}q^{\binom{j}{2}}{m\brack j}.
\end{equation}
\end{lem}

\pf Appealing to Lemma \ref{rr-1-new}, we get
\begin{align*}
\sum_{\pi\in\mathcal{C}_2}q^{|\pi|}&=1+\sum_{m\geq1}\frac{1}{(q;q)_{m}}\sum_{j\geq0}\left[b_\mathcal{G}(m,{m+j},{1})+b_\mathcal{G}(m,m+j-1,\bar{1})\right]\\[5pt]
&=\sum_{m\geq0}\frac{1}{(q;q)_{m}}
\sum_{j\geq0}\left(q^{\binom{m+1}{2}+\binom{j+1}{2}}{{m-1}\brack j}+q^{\binom{m+1}{2}+\binom{j}{2}}{{m-1}\brack {j-1}}\right)\\[5pt]
&=\sum_{m\geq0}\frac{q^{\binom{m+1}{2}}}{(q;q)_{m}}
\sum_{j\geq0}q^{\binom{j}{2}}\left(q^j{{m-1}\brack j}+{{m-1}\brack {j-1}}\right)\\[5pt]
&=\sum_{m\geq0}\frac{q^{\binom{m+1}{2}}}{(q;q)_{m}}
\sum_{j\geq0}q^{\binom{j}{2}}{{m}\brack j},
\end{align*}
where the final equation follows from \eqref{bin-r-1}. Thus, we complete the proof.  \qed

We are now in a position to give a proof of \eqref{main-2}.

{\noindent \bf Proof of \eqref{main-2}.} Summing up  the $j$-sum in the right-hand side of \eqref{rr-1-main} by letting $z=-1$ and $N=m$ in \eqref{j-sum}, we get
\begin{equation}\label{rr-1-j-sum-1}
\sum_{\pi\in\mathcal{C}_2}q^{|\pi|}=\sum_{m\geq0}\frac{(-1;q)_mq^{\binom{m+1}{2}}}{(q;q)_{m}}.
\end{equation}

On the other hand, by interchanging the order of summation in the right-hand side of \eqref{rr-1-main}, we can obtain that
\begin{align*}\label{rr-1-j-sum-2}
\sum_{\pi\in\mathcal{C}_2}q^{|\pi|}&=\sum_{j\geq0}\frac{q^{\binom{j}{2}}}{(q;q)_j}\sum_{m\geq j}\frac{q^{\binom{m+1}{2}}}{(q;q)_{m-j}}\\
&=\sum_{j\geq0}\frac{q^{\binom{j}{2}}}{(q;q)_j}\sum_{m\geq 0}\frac{q^{\binom{m+j+1}{2}}}{(q;q)_{m}}\\
&=\sum_{j\geq0}\frac{q^{j^2}}{(q;q)_j}\sum_{m\geq 0}\frac{q^{\binom{m}{2}+m(j+1)}}{(q;q)_{m}}\\
&=\sum_{j\geq0}\frac{q^{j^2}}{(q;q)_j}(-q^{j+1};q)_\infty\qquad\text{by \eqref{Euler-2}}\\
&=(-q;q)_\infty\sum_{j\geq0}\frac{q^{j^2}}{(q^2;q^2)_j}\\
&=\frac{(-q;q^2)_\infty}{(q;q^2)_\infty}\qquad\text{by \eqref{Euler-2}}.
\end{align*}

Combining with \eqref{rr-1-j-sum-1}, we arrive at \eqref{main-2}. This completes the proof.  \qed

\subsection{Overpartition analogue of the second Rogers-Ramanujan identity}

In this subsection, we focus on $\mathcal{G}_1(m)$ and  the proof of \eqref{main-3}. For $m\geq1$ and $j\geq1$, let $b_\mathcal{G}(m,j,2)$ be the generating function for the overpartitions of $\mathcal{G}(m)$ with the largest part of size $j$ and the smallest part $2$. Clearly,
\[\sum_{\pi\in\mathcal{G}_1(m)}q^{|\pi|}=\sum_{j\geq 1}\left[b_\mathcal{G}(m,j,\bar{1})+b_\mathcal{G}(m,j,2)\right].\]

\begin{lem}\label{rr-2-new}
For $m\geq 1$ and $j\geq 0,$
\begin{equation}\label{rr-2-eqn}
b_\mathcal{G}(m,m+j,2)=q^{m+\binom{m+1}{2}+\binom{j}{2}}
{{m-1}\brack {j-1}}.
\end{equation}
\end{lem}

\pf For $m\geq 1$ and $j=0$, assume that $\pi=(\pi_1,\pi_2,\ldots,\pi_m)$ is an overpartition in $\mathcal{G}(m)$ with the largest part of size $m+j$ and the smallest part $2$, that is, $\pi_1=m+j$ or $\overline{m+j}$ and $\pi_m=2$. Then, we have
\[\pi_1\geq \pi_2+1\geq\cdots\geq \pi_m+(m-1)=m+1.\]
This implies that $j\geq1$, and so $b_\mathcal{G}(m,m,2)=0$. Hence, \eqref{rr-2-eqn} holds for $m\geq 1$ and $j=0$.

For $m\geq 1$ and $j\geq 1$, we need to show that
\begin{equation}\label{lem-rr-2-pr}
b_\mathcal{G}(m,m+j,2)=q^mb_\mathcal{G}(m,m+j-1,1)
\end{equation}

If we subtract $1$ from each part of $\pi$, then we get an overpartition in $\mathcal{G}(m)$ with the largest part of size $m+j-1$ and the smallest part $1$, and vice {versa}. So, \eqref{lem-rr-2-pr} is valid. Combining \eqref{rr-2-pr} and \eqref{lem-rr-2-pr}, we arrive at
\eqref{rr-2-eqn}. This completes the proof.   \qed

Now, we give the generating function for the overpartitions in $\mathcal{C}_1$.
\begin{lem}
\begin{equation}\label{rr-2-main-0}
\sum_{\pi\in\mathcal{C}_1}q^{|\pi|}=\sum_{m\geq0}\frac{q^{\binom{m+1}{2}}}{(q;q)_{m}}
\sum_{j\geq0}q^{\binom{j+1}{2}}{m\brack j}.
\end{equation}
\end{lem}

\pf By Lemmas \ref{rr-1-new} and \ref{rr-2-new}, we get
\begin{align*}
\sum_{\pi\in\mathcal{C}_1}q^{|\pi|}&=1+\sum_{m\geq1}\frac{1}{(q;q)_{m}}\sum_{j\geq0}\left[b_\mathcal{G}(m,{m+j},\bar{1})+b_\mathcal{G}(m,m+j,2)\right]\\[5pt]
&=1+\sum_{m\geq1}\frac{1}{(q;q)_{m}}
\sum_{j\geq0}\left(q^{\binom{m+1}{2}+\binom{j+1}{2}}
{{m-1}\brack j}+q^{m+\binom{m+1}{2}+\binom{j}{2}}
{{m-1}\brack {j-1}}\right)\\[5pt]
&=1+\sum_{m\geq1}\frac{q^{\binom{m+1}{2}}}{(q;q)_{m}}
\sum_{j\geq0}q^{\binom{j+1}{2}}\left({{m-1}\brack j}+q^{m-j}{{m-1}\brack {j-1}}\right)\\[5pt]
&=\sum_{m\geq0}\frac{q^{\binom{m+1}{2}}}{(q;q)_{m}}
\sum_{j\geq0}q^{\binom{j+1}{2}}{m\brack j},
\end{align*}
where the final equation follows from \eqref{bin-r-2}. Thus, we complete the proof.  \qed

We are now in a position to give a proof of \eqref{main-3}.

{\noindent \bf Proof of \eqref{main-3}.} Summing up  the $j$-sum in the right-hand side of \eqref{rr-1-main} by letting $z=-q$ and $N=m$ in \eqref{j-sum}, we get
\begin{equation}\label{rr-2-j-sum-1-0}
\sum_{\pi\in\mathcal{C}_1}q^{|\pi|}=\sum_{m\geq0}\frac{(-q;q)_mq^{\binom{m+1}{2}}}{(q;q)_{m}}.
\end{equation}

On the other hand, by interchanging the order of summation in the right-hand side of \eqref{rr-2-main-0}, we can obtain that
\begin{align*}\label{rr-2-j-sum-2-0}
\sum_{\pi\in\mathcal{C}_1}q^{|\pi|}&=\sum_{j\geq0}\frac{q^{\binom{j+1}{2}}}{(q;q)_j}\sum_{m\geq j}\frac{q^{\binom{m+1}{2}}}{(q;q)_{m-j}}\\[5pt]
&=\sum_{j\geq0}\frac{q^{\binom{j+1}{2}}}{(q;q)_j}\sum_{m\geq 0}\frac{q^{\binom{m+j+1}{2}}}{(q;q)_{m}}\\[5pt]
&=\sum_{j\geq0}\frac{q^{j^2+j}}{(q;q)_j}\sum_{m\geq 0}\frac{q^{\binom{m}{2}+m(j+1)}}{(q;q)_{m}}\\[5pt]
&=\sum_{j\geq0}\frac{q^{j^2+j}}{(q;q)_j}(-q^{j+1};q)_\infty\qquad\text{by \eqref{Euler-2}}\\[5pt]
&=(-q;q)_\infty\sum_{j\geq0}\frac{q^{j^2+j}}{(q^2;q^2)_j}\\[5pt]
&=(-q;q)_\infty(-q^2;q^2)_\infty\qquad\text{by \eqref{Euler-2}}\\[5pt]
&=\frac{(-q;q)_\infty}{(q^2;q^4)_\infty}.
\end{align*}

Combining with \eqref{rr-2-j-sum-1-0}, we arrive at \eqref{main-3}. This completes the proof.  \qed

{F \noindent{\bf Acknowledgments.}
This work was supported by
the National Science Foundation of
China  (Nos. 12101437).  We are greatly indebted to the referees for their insightful suggestions leading to
an improvement of an earlier version.}


\begin{thebibliography}{99}

 \bibitem{Andrews-1976} G.E. Andrews, The Theory of partitions, Addison-Wesley Publishing Co., 1976.

 \bibitem{Andrews-2018} G.E. Andrews, Integer partitions with even parts below odd parts and the mock theta functions, Ann. Comb. 22 (2018) 433--445.

\bibitem{Andrews-2019} G.E. Andrews, Partitions with parts separated by parity, Ann. Comb. 23 (2019) 241--248.

 \bibitem{Andrews-2022} G.E. Andrews, Separable integer partition classes, Trans. Amer. Math. Soc. 9 (2022) 619--647.


 \bibitem{Chen-Sang-Shi-2013} W.Y.C. Chen, D.D.M. Sang and D.Y.H. Shi, The Rogers-Ramanujan-Gordon theorem for overpartitions, Proc. London Math. Soc. 106  (2013) 1371--1393.



     \bibitem{Corteel-Lovejoy-2004} S. Corteel and J. Lovejoy, Overpartitions, Trans. Amer. Math. Soc. 356 (4) (2004) 1623--1635.


\bibitem{Lovejoy-2003} J. Lovejoy, Gordon's theorem for overpartitions, J. Combin. Theory, Ser. A 103 (2003) 393--401.

\bibitem{Passary-2019} D. Passary, Studies of partition functions with conditions on parts and parity, PhD thesis, Penn. State University, 2019.

\end{thebibliography}
\end{document}